\documentclass[letterpaper,12pt]{article}
\usepackage{tabularx} % extra features for tabular environment
\usepackage{amsmath}  % improve math presentation
\usepackage{graphicx} % takes care of graphic including machinery
\usepackage{amssymb}
\usepackage[margin=1in,letterpaper]{geometry} % decreases margins
\usepackage{cite} % takes care of citations
\usepackage[final]{hyperref} % adds hyper links inside the generated pdf file
\hypersetup{
	colorlinks=true,       % false: boxed links; true: colored links
	linkcolor=black,        % color of internal links
	citecolor=black,        % color of links to bibliography
	filecolor=magenta,     % color of file links
	urlcolor=blue         
}
\begin{document}

\title{\textbf{Weak inverse problem of calculus of variations for geodesic mappings and relation to harmonic maps}}
\author{Stanislav Hronek}
\date{}

\maketitle
\textit{Department of Theoretical Physics and Astrophysics, Masaryk University, 611 37 Brno, Czech Republic}

\begin{abstract}
In this paper we study the relation between geodesic and harmonic mappings. Harmonic mappings are defined between Riemannian manifolds as critical points of the energy functional, on the other hand geodesic mappings are defined in a more general setting (manifolds with affine connections). Using the well-established formalism of calculus of variations on fibred manifolds we solve the weak inverse problem for the equation of geodesic mappings and get a variational equation which is a consequence of the geodesic mappings equation. For the connection on the target manifold we get the expected result, that it is a metric connection. However we find that the connection on the source manifold need not be metric. The interesting result is that the metric which induces the connection on the target manifold can change between fibres and these changes are related to the connection on the source manifold. These results hint onto a possibility for a more general structure on the fibred manifold, than usually assumed.
\end{abstract}

\section{Geodesic mappings and basic setting} 
\label{sec:1}
Let us start with geodesic mappings of manifolds with affine connections. For the theory of geodesic mappings we refer to \cite{mik}. Because geodesics on manifolds are characterized by the symmetric part of the connection only, we can restrict ourselves to torsion-free manifolds with affine connections, i.e from now on, we assume that all connections under consideration are symmetric. We will also only consider naturally parametrized geodesics. Consider manifolds with affine connections  $(M,\,^M\nabla)$ and $(N,\,^N\nabla)$ and a map between them $\phi: (M,\,^M\nabla)\longrightarrow(N,\,^N\nabla)$. This map is said to be a geodesic map if
\begin{enumerate}
    \item  $\phi$ is a diffeomorphism of $M$ onto $N$; and
    \item the image under $\phi$ of any geodesic arc in $M$ is a geodesic arc in $N$; and
    \item the image under the inverse function $\phi^{-1}$ of any geodesic arc in $N$ is a geodesic arc in $M$.
\end{enumerate}
The usual example of a geodesic mapping would be an isometry of euclidean surfaces.
In our paper we will generalize this definition a little, we will give up the assumptions $1.$ and $3.$ meaning instead of diffeomorphisms we will be working with immersion. Mathematically speaking a mapping $\phi: (M,\,^M\nabla)\longrightarrow(N,\,^N\nabla)$, $(dim (M)\leq dim (N))$ is geodesic if for every geodesic curve $x(t)$ on $(M,\,^M\nabla)$, $\phi \circ x(t)$ is a geodesic curve on $(N,\,^N\nabla)$. For solving the inverse problem of calculus of variation we would like to use the formalism of calculus of variations on fibred manifolds. The fibred space for the problem will be $(M\times N,\pi, M)$, which has dimension $m+n$ and the mapping $\phi$ now serves as a fiber coordinate. We also suppose geodesics on $M$ and $N$ are parametrized by the same parameter $t$. From the simple definition of a geodesic mappings one can derive a set of geodesic equations which will serve as conditions for the mapping $\phi$ to be a geodesic mapping. Let us write out the equations in coordinate systems $(x^i, \phi^\sigma)$ on the total space and adapted system on the basis $(x^i)$. In these coordinate systems the affine connections $^M\nabla$ respectively $^N\nabla$ have components denoted by $^M\Gamma^h_{ij}$ respectively $^N\Gamma^\mu_{\nu\lambda}$. The geodesic equations for a curve $x(t)$ in $M$ and a curve $y(t)$ in $N$ are
\begin{align*}
\ddot{x}^h+\,^M\Gamma^h_{ij}\dot{x}^i\dot{x}^j&=0,\\
\ddot{y}^\sigma+\,^N\Gamma^\sigma_{\mu\nu}\dot{y}^\mu\dot{y}^\nu&=0,
\end{align*}
$$
i,j,l=1,\ldots,m=dim(M) \quad \mu,\nu,\lambda=1,\ldots,n=dim(N),
$$
For a geodesic mapping $\phi$ the geodesic curve $y(t)$ is the image of $x(t)$ by $\phi$, substituting  $y(t)=\phi(x(t))$ in the second equation we get
\begin{align*}
\frac{{\rm d}}{{\rm d}t}\left(\phi^\mu_l\dot{x}^l\right)+^N\Gamma^\mu_{\nu\lambda}\phi^\nu_i\dot{x}^i\phi^\lambda_j\dot{x}^j&=0,
\end{align*}
where we used the chain rule in the second equation $\frac{{\rm d}}{{\rm d}t}(\phi^\mu(x^l(t)))=\phi^\mu_l\dot{x}^l$. Computing the derivative in the second equation and then substituting for the second derivative from the first we get
$$
\frac{{\rm d}}{{\rm d}t}\left(\phi^\mu_l\dot{x}^l\right)+^N\Gamma^\mu_{\nu\lambda}\phi^\nu_i\dot{x}^i\phi^\lambda_j\dot{x}^j=\phi^\mu_{kl}\dot{x}^k\dot{x}^l+\phi^\mu_l\ddot{x}^l+^N\Gamma^\mu_{\nu\lambda}\phi^\nu_i\dot{x}^i\phi^\lambda_j\dot{x}^j=0,
$$
$$
\phi^\mu_{kl}\dot{x}^k\dot{x}^l - ^M\Gamma^h_{ij}\dot{x}^i\dot{x}^j\phi^\mu_h+^N\Gamma^\mu_{\nu\lambda}\phi^\nu_i\dot{x}^i\phi^\lambda_j\dot{x}^j=\dot{x}^i\dot{x}^j\left(\phi^\mu_{ij}-^M\Gamma^h_{ij}\phi^\mu_h+^N\Gamma^\mu_{\nu\lambda}\phi^\nu_i\phi^\lambda_j\right),
$$
\begin{equation}
\label{rce}
\phi^\sigma_{ij}-^M\Gamma^k_{ij}\phi^\sigma_k+^N\Gamma^\sigma_{\alpha\lambda}\phi^\alpha_i\phi^\lambda_j=0.
\end{equation}
What we get is a sufficient condition for $\phi$ to be a geodesic mapping. The second part of interest is harmonic mappings.

\section{Harmonic mappings}
\label{sec:2}
The basics of harmonic mappings can be found in \cite{eel} and \cite{ura}. Their applications to physics, which include string theory, sigma models and general relativity, are presented in papers \cite{mis} and \cite{nor}. The main change from geodesic mappings is in the setting. Harmonic mappings are defined on Riemannian manifolds, i.e manifolds endowed with a metric tensor. We say that a mapping $\phi$ between two Riemannian manifolds $(M,g)$ and $(N,h)$ is harmonic if it is a stationary (extremal) point of the energy functional.
$$
E(\phi)=\int_M \frac{1}{2}Tr_g(\phi^* h) \omega_0,
$$
where $\omega_0$ is the volume element on $M$ corresponding to the metric tensor $g$. Euler-Lagrange equations of this functional yield similar equations as for geodesic mappings (\ref{rce}), with the difference that for harmonic mappings the connections are metric connections.
The Lagrange function in the chosen coordinate system $(x^i,\phi^\sigma)$ has the following form
$$
L=\frac{1}{2}g^{ij}h_{\alpha\lambda}\phi^\alpha_i\phi^\lambda_j
$$
Its Euler-Lagrange equations are as follows
$$
{\rm d}_k\frac{\partial L}{\partial \phi^\sigma_k}=\frac{\partial L}{\partial \phi^\sigma}
$$
After computing the derivatives and arranging the terms we get
$$
g^{ij}h_{\sigma\nu}\phi^\sigma_{ij}+g^{ij}\phi^\alpha_i\phi^\lambda_j\left(\frac{1}{2}h_{\sigma\lambda,\alpha}+\frac{1}{2}h_{\sigma\alpha,\lambda}-\frac{1}{2}h_{\alpha\lambda,\sigma}\right)+g^{kj}_{\,\,\,\,,k}h_{\sigma\nu}\phi^\sigma_j=0
$$
%relabeling sigma->nu
\begin{equation}
\label{har}
g^{ij}h_{\sigma\nu}\left(\phi^\sigma_{ij}-^M\Gamma^k_{ij}\phi^\sigma_k+^N\Gamma^\sigma_{\alpha\lambda}\phi^\alpha_i\phi^\lambda_j\right)=0,
\end{equation}
where in the last step we used an expression for the metric trace of Christoffel symbols $g^{ij}\Gamma^k_{ij}=-g^{kl}_{\,\,\,\,,l}$. These equations seem rather complicated, let us give some examples of harmonic mappings. Constant maps are harmonic. If the source manifold $(M,g)$ was
$\mathbb{R}$ with the natural metric, the mapping $\phi$ would be harmonic if and only if it was a geodesic. On the other hand if $\mathbb{R}$ with the natural metric was the target space, the equation (\ref{har}) would be the Laplace equation and its solutions, called harmonic functions, are a special case of harmonic mappings. As mentioned above we get similar equations, more concretely we get a trace of the geodesic equations (\ref{rce}). There is also additional difference being affine connections on one hand and metric on the other. The similarities suggest the following question: what is the connection between variationality of equations (\ref{rce}) and the corresponding connections being metric ?

% Always give a unique label
% and use \ref{<label>} for cross-references
% and \cite{<label>} for bibliographic references
% use \sectionmark{}
% to alter or adjust the section heading in the running head

\section{Weak inverse problem of calculus of variations}
\label{sec:3}
To answer this question we will start with the equations (\ref{rce}) assume that the connections are affine and not necessary metric and solve the inverse problem of calculus of variations. From the equations it is obvious that they are not variational by themselves, meaning we need to impose and solve the weak inverse problem. To do that we multiply the equation by a multiplier $B^{ij}_{\sigma\nu}(x^k,\phi^\mu)$ and assume it does not depend on the derivatives of the basis and fibre coordinates, which is usually the case, but this also means our conclusions will be only sufficient and possibly not-necessary. By multiplying the equations (\ref{rce}) by a multiplier $B$ we get new equations which do not need to be equivalent to (\ref{rce}) but are a differential consequence of (\ref{rce}). 

We associate a dynamical $(m+1)$-form with the equation, this form is the following $E=E_\nu\omega^\nu\wedge\omega_0$, where $\omega_0$ is the volume element on $M$ and $\omega^\nu={\rm d}\phi^\nu-\phi^\nu_i{\rm d}x^i$ is the contact 1-form on $M\times N$.
\begin{equation}
\label{EL}
E_\nu=B_{\sigma\nu}^{ij}\left(\phi^\sigma_{ij}-^M\Gamma^k_{ij}\phi^\sigma_k+^N\Gamma^\sigma_{\alpha\lambda}\phi^\alpha_i\phi^\lambda_j\right).
\end{equation}

Again using the logic that we are studying geodesic equations and the connections and their components are supposed to be symmetric we can assume the multiplier is also symmetric in the upper indices $i, j$.

Instead of solving variationality the equation (\ref{rce}) we study the variationality of the associated form $E$ using the tools of calculus of variations on fibred manifolds, which can be found in \cite{dem} and \cite{kru}. The conditions for this form to be variational are called Helmholtz conditions of variationality, which are of the following form (derivation can be found in \cite{kru})
\begin{align}
\label{1HP}
\frac{\partial E_\nu}{\partial \phi^\mu_{lp}}-\frac{\partial E_\mu}{\partial \phi^\nu_{lp}}&=0,\\
\label{2HP}
\frac{\partial E_\nu}{\partial \phi^\mu_{l}}+\frac{\partial E_\mu}{\partial \phi^\nu_{l}}-2{\rm d}_p\frac{\partial E_\mu}{\partial \phi^\nu_{lp}}&=0,\\
\label{3HP}
\frac{\partial E_\nu}{\partial \phi^\mu}-\frac{\partial E_\mu}{\partial \phi^\nu}+{\rm d}_l\frac{\partial E_\mu}{\partial \phi^\nu_{l}}-{\rm d}_l{\rm d}_p\frac{\partial E_\mu}{\partial \phi^\nu_{lp}}&=0.
\end{align}

From the first condition we get symmetry of the multiplier B in the lower indices
$$
B^{ij}_{\sigma\nu}=B^{ij}_{\nu\sigma}.
$$
Because at the beginning we assumed B does not depend on derivatives the equations for the second condition split into a polynomial form in the derivatives of $\phi$. Setting each coefficient to zero results in two conditions.
\begin{align}
\label{HP21}
-B^{ij}_{\mu\nu}\,^M\Gamma^l_{ij}&=\frac{\partial B^{lp}_{\mu\nu}}{\partial x^p}.
\\
\label{HP22}
^N\Gamma^\sigma_{\mu\lambda}B^{ij}_{\sigma\nu}+\,^N\Gamma^\sigma_{\nu\lambda}B^{ij}_{\sigma\mu}&=\frac{\partial B^{ij}_{\mu\nu}}{\partial \phi^\lambda},
\end{align}
These conditions already tell us something about the form of the multiplier B. The second equation is the condition for a connection $^N\Gamma^\sigma_{\mu\lambda}$ to be to be compatible with metric tensors $B^{ij}$, with components $B^{ij}_{\sigma\nu}$ for any choice of indices $i,j$. Noticing that the equations separate, in the sense that in the equation (\ref{HP21}) there is a derivative with respect to $x^p$ and in the equation (\ref{HP22}) with respect to $\phi^\lambda$, we can guess that the multiplier $B$ separates also into the following form
$$
B^{ij}_{\sigma\nu}=g^{ij}(x^k)h_{\sigma\nu}(\phi^\mu),
$$
This would be the simplest choice, in particular we know this choice is correct (is a solution to the inverse problem) because it gives the harmonic mappings equation (\ref{har}). We can be more general and allow the functions $h_{\sigma\nu}(\phi^\mu)$ to also depend on $x^k$. The reasoning is the following. The second equation tells us that $h_{\sigma\nu}$ are components of a metric tensor of the connection $^N\Gamma^\sigma_{\mu\lambda}$ at one particular fibre $\pi^{-1}(x^k)$. If we move to another fibre the connection $^N\Gamma^\sigma_{\mu\lambda}$ is also metric but the metric can be different. We can have a different metric tensor $h$ in each fibre and we allow the functions $h_{\sigma\nu}$ to also depend on $x^k$ for this very reason.  
Therefore we choose the multiplier in the following form
\begin{equation}
\label{mult}
B^{ij}_{\sigma\nu}=g^{ij}(x^k)h_{\sigma\nu}(x^k,\phi^\mu).
\end{equation}
To justify calling the functions $g^{ij}$ and $h_{\sigma\nu}$ components of metric tensors we need to check if they are symmetric and regular. Their symmetry follows from symmetry of the multiplier $B$, symmetry of $B$ in upper indices we assumed and in lower indices we got from Helmholtz conditions, the regularity also follows from regularity of the multiplier $B$. In coordinate-free form we have 
$
B=g\otimes h.
$

We know that the equation (\ref{HP22}) assures that the connection $^N\Gamma^\sigma_{\mu\lambda}$ comes from a metric $h$, which can be different in each fibre (for different $x$). We can calculate how it changes between fibres from the equation (\ref{HP21}). We substitute for $B$ into the equation (\ref{HP21}) and simplify
\begin{align*}
&-B^{ij}_{\mu\nu}\,^M\Gamma^l_{ij}=\frac{\partial B^{lp}_{\mu\nu}}{\partial x^p}\\
&-g^{ij}h_{\mu\nu}\,^M\Gamma^l_{ij}=g^{lp}_{\,\,\,,p}h_{\mu\nu}+h_{\mu\nu,p}g^{lp}\\
&h_{\mu\nu}\left(g^{ij}\,^M\Gamma^l_{ij}+g^{lp}_{\,\,\,,p}\right)=-h_{\mu\nu,p}g^{lp}
\end{align*}
where $g^{lp}_{\,\,\,,p}$ is actually a trace of a connection induced by the metric tensor $g$ let us denote it by $^M\bar{\nabla}$, meaning the equation is a difference of traces of two connections. Originally we assumed the space $M$ is only endowed with an affine connection, however we see there also supposedly exists a metric $g$, but that is no surprise because we know every smooth manifold admits a metric. We see that the dependency of metric tensor $h$ on the basis coordinate $x$ is given by the difference of traces of connections on the space $M$.
We can then express the relation between the connection components 
$$
^M\Gamma^l_{ij}= \, ^M\bar{\Gamma}^l_{ij}+S^l_{ij},
$$
where $S^l_{ij}$ is a tensor which satisfies
$$
g^{ij}S^l_{ij}=-\frac{1}{n}h^{\mu\nu}h_{\mu\nu,p}g^{lp}.
$$

The last remaining Helmholtz condition unfortunately bring no new information after rearranging it into a polynomial form and from requiring that all the polynomial coefficients vanish we get four conditions 
\begin{align}
\label{HP31}
&\partial_\nu B^{ij}_{\sigma\mu}-\partial_\mu B^{ij}_{\sigma\nu}-\partial_\sigma B^{ij}_{\mu\nu}+B^{ij}_{\alpha\nu}\, ^N\Gamma^\alpha_{\mu\sigma}+B^{ij}_{\alpha\nu}\, ^N\Gamma^\alpha_{\sigma\mu}=0\\
\nonumber
&\partial_\nu (B^{ij}_{\sigma\mu})^N\Gamma^\sigma_{\alpha\lambda}+B^{ij}_{\sigma\mu}\partial_\nu^N\Gamma^\sigma_{\alpha\lambda}-\partial_\mu(B^{ij}_{\sigma\nu})^N\Gamma^\sigma_{\alpha\lambda}-B^{ij}_{\sigma\nu}\partial_\mu^N\Gamma^\sigma_{\alpha\lambda}+\partial_\alpha(B^{ij}_{\sigma\nu})^N\Gamma^\sigma_{\mu\lambda}+B^{ij}_{\sigma\nu}\partial_\alpha^N\Gamma^\sigma_{\mu\lambda}\\
\label{HP32}
&+\partial_\lambda(B^{ij}_{\sigma\nu})^N\Gamma^\sigma_{\alpha\mu}+B^{ij}_{\sigma\nu}\partial_\lambda\,^N\Gamma^\sigma_{\alpha\mu}-\partial_\lambda\partial_\alpha B^{ij}_{\mu\nu}=0\\
\label{HP33}
&^M\Gamma^k_{ij}\left(\partial_\mu B^{ij}_{\sigma\nu}-\partial_\nu B^{ij}_{\sigma\mu}-\partial_\sigma B^{ij}_{\mu\nu}\right)
+2\partial_l\left(B^{lk}_{\lambda\nu}\right)\,^N\Gamma_{\mu\sigma}^\lambda-2\partial_p\partial_\sigma B^{kp}_{\mu\nu}=0
\\
\label{HP34}
&-\partial_l(B^{ij}_{\mu\nu})^M\Gamma^l_{ij}-B^{ij}_{\mu\nu}\partial_l^M\Gamma^l_{ij}-\partial_l\partial_pB^{lp}_{\mu\nu}=0.
\end{align}
The equation (\ref{HP31}) provides us with the same information as (\ref{HP22}) that being, the connections $^N\nabla$ is metric with the metric tensor $h$.

We can also see that the remaining equations are dependent. The equation (\ref{HP34}) is just a derivative of (\ref{HP21}) with respect to $x^l$ and equation (\ref{HP33}) 
is a multiple of the equation (\ref{HP31}) by $^M\Gamma^k_{ij}$. The only equation that does not depend on the previous ones is the equation (\ref{HP32}) but after some calculations it results in an identity for the Riemann curvature tensor $R$. 

$$
h_{\beta\mu}\,R^\beta_{\,\,\lambda\nu\alpha}=h_{\sigma\nu}\,R^\sigma_{\,\,\alpha\mu\lambda}\longrightarrow R_{\mu\lambda\nu\alpha}=R_{\nu\alpha\mu\lambda}.
$$

\section{Summary and conclusions}
\label{sec:4}

We now summarize everything we discovered from the Helmholtz conditions.
We have the equation (\ref{rce})
for a geodesic mapping $\phi$. We ask the following question: what is the connection between variationality of this equation and the corresponding connections being metric ?  Therefore we are solving an inverse problem for the associated dynamical form $E=E_\nu\, \omega^\nu\wedge\omega_0$
$$
E_\nu=B^{ij}_{\sigma\nu}\left(\phi^{\sigma}_{ij}- ^M \Gamma^k_{ij}\phi^\sigma_k+^N\Gamma^\sigma_{\alpha\lambda}\phi^\alpha_i\phi^\lambda_j\right)
$$
We choose a specific form (\ref{mult}) of the variational multiplier $B$. The conditions for variationality are
\begin{enumerate}
\item Connection $^N\nabla$ is metric and is fiber-wise induced by the metric $h$. (Metric $h$ is generally different in each fibre $h_{\sigma\nu}=h_{\sigma\nu}(x^k,\phi^\mu)$. The way in which this metric changes in $x^k$ is given by the connection $^M \nabla$ and metric $g$)
\item Connection $^M \nabla$ does not need to be metric but is related to a metric connection by
$$
\Gamma^k_{ij}=\bar{\Gamma}^k_{ij}+S^k_{ij},
$$
where $S^k_{ij}$ is a tensor whose metric trace by the tensor $g$ relates to the changes of the metric $h$  
$$
g^{ij}S^l_{ij}=-h^{\mu\nu}h_{\mu\nu,p}g^{lp}.
$$
\end{enumerate}
The form (\ref{mult}) of the variational multiplier $B$ is a solution to the inverse problem if the corresponding metrics satisfy the above conditions. The results suggests that both spaces are Riemannian, where $h$ is compatible with the connection $^N\nabla$ on $N$ but $g$ is not necessarily compatible with $^M \nabla$ on $M$. The interesting result is that the metric $h$ can be different in different fibres and the changes are related to structures on the base manifold $M$. This conclusion mostly results from using the formalism of fibred manifolds and finding a non-trivial solution to the weak inverse problem, but it suggests more complicated fibred structure for the problem of geodesic and harmonic mappings. There remains the question of finding more general forms of the multiplier $B$, which could be a part of further research.

\end{document}